\documentclass[10pt]{amsart}
\usepackage{amssymb}
\numberwithin{equation}{section}

\newtheorem{theorem}{Theorem}[section]

\newtheorem{corollary}[theorem]{Corollary}
\newtheorem{proposition}[theorem]{Proposition}
\newtheorem{lemma}[theorem]{Lemma}
\theoremstyle{definition}
\newtheorem{definition}[theorem]{Definition}
\theoremstyle{remark}
\newtheorem{remark}{Remark}
\newtheorem{example}{Example}

\newcommand{\cV}{\mathcal{V}}
\newcommand\cH{\mathcal{H}}
\newcommand\cP{\mathcal{P}}
\newcommand\cR{\mathcal{R}}
\newcommand\cS {\mathcal{S}}

\newcommand{\R}{\mathbb{R}}
\newcommand{\C}{\mathbb{C}}
\newcommand{\Z}{\mathbb{Z}}

\newcommand{\bbT}{\mathbb{T}}

\newcommand\lie[1]{\mathfrak{#1}}

\newcommand{\fg}{\lie{g}}
\newcommand{\fm}{\lie{m}}

\newcommand{\ft}{\lie{t}}

\def	\inv	{^{-1}}

\newcommand\Diff{\mathop{\it Diff}\nolimits}

\newcommand\Curv{\mathop{\it Curv }\nolimits}
\newcommand\curv{\mathop{\it Curv }\nolimits}

\begin{document}

\title{Contact fiber bundles}

\author{Eugene Lerman}
\address{Department of
Mathematics, University of Illinois, Urbana, IL 61801}
\email{lerman@math.uiuc.edu}

\thanks{Supported by in part by the Swiss NSF, US NSF grant DMS-0204448 and by 
R. Kantorovitz.}
\date{\today}

\begin{abstract}
We define contact fiber bundles and investigate conditions for the
existence of contact structures on the total space of such a  bundle.
The results are analogous to minimal coupling in symplectic geometry.
The two applications are  construction of $K$-contact manifolds
generalizing Yamazaki's fiber join construction and a cross-section
theorem for contact moment maps.
\end{abstract}

\maketitle

\tableofcontents

\section{Introduction }

A few years ago I came across an interesting paper by Yamazaki
\cite{Yamazaki} in which $K$-contact manifolds were constructed ``by
fiber join.''  The contact manifolds in question were odd-dimensional
sphere bundles over Riemann surfaces.  These bundles were associated
to certain principal torus bundles, and the choice of a contact form
on the total space involved a choice of a connection on the torus
bundle whose curvature had to satisfy certain non-degeneracy
condition.  It all very much looked like a contact version of
Sternberg's minimal coupling construction in symplectic geometry.  The
goal of the present paper is to explain why this is indeed the case.
In order to do this systematically I felt it is necessary to first sort out
the definition of a contact fiber bundle and  to investigate
conditions for the existence of contact structures on the total space
of such a bundle.  Of course, since any co-oriented contact manifold
is the quotient of a symplectic cone by dilations and since any
symplectic cone is a symplectization of a contact manifold, one can
argue that contact fiber bundles are simply symplectic fiber bundles
with fibers being symplectic cones.  However, it seems useful to study
the matter completely in contact terms.  For example, the notion of a
$K$-contact structure does not translate naturally into symplectic
terms.

\begin{remark}\label{remark1}
In this paper we assume that all contact structures are co-oriented.
Recall that a contact structure $\xi $ on a manifold $M$ is
{\bf co-oriented} if there exists a 1-form $\alpha$ with $\ker \alpha =
\xi$.  Equivalently $\xi$ is co-oriented if its annihilator $\xi
^\circ \subset T^*M$ is an oriented line bundle. That is, the line
bundle minus the zero section, $\xi ^\circ \smallsetminus M$, has two
components and we have single out one of the components, call it
$\xi^\circ_+$.  We thus may think of $\xi^\circ_+$ as a co-orientation
of $\xi$. Note that the image of $\alpha : M \to TM$ with $\ker \alpha =
\xi$ singles out one of the components of $\xi^\circ \smallsetminus
M$.  The same remark applies to any codimension 1 distribution $\xi$
on $M$, not just to contact ones.  That is, we only consider
co-oriented distributions.
\end{remark}

\subsection*{A note on notation}  

If $U$ is a subspace of a vector space $V$ we denote its annihilator
in the dual vector space $V^*$ by $U^\circ$.  Thus $U^\circ = \{ \ell
\in V^* \mid \ell |_U = 0\}$.  We use the same notation for distributions.

Throughout the paper the Lie algebra of a Lie group denoted by a
capital Roman letter is denoted by the same small letter in the
fraktur font: thus $\fg$ denotes the Lie algebra of a Lie group $G$
etc.  The vector space dual to $\fg$ is denoted by $\fg^*$. The
identity element of a Lie group is denoted by 1.  The natural pairing
between $\fg$ and $\fg^*$ is denoted by $\langle \cdot, \cdot
\rangle$.

When a Lie group $G$ acts on a manifold $M$ we denote the action by an
element $g\in G$ on a point $x\in M$ by $g\cdot x$; $G\cdot x$ denotes
the $G$-orbit of $x$ and so on.  The vector field induced on $M$ by an
element $X$ of the Lie algebra $\fg$ of $G$ is denoted by $X_M$ 
(that is, $X_M (x) = \left. \frac{d}{dt} \right|_0 (\exp tX) \cdot x$)
and the diffeomorphism induced by $g\in G$ on $M$ by $g_M$.  Thus in this
notation $g\cdot x = g_M (x)$.  The isotropy group of a point $x\in M$
is denoted by $G_x$; the Lie algebra of $G_x$ is denoted by $\fg_x$
and is referred to as the isotropy Lie algebra of $x$.  Recall that
$\fg_x = \{ X \in \fg\mid X_M (x) = 0\}$.

If $X$ is a vector field and $\tau$ is a tensor, then $L_X \tau$
denotes the Lie derivative of $\tau$ with respect to $X$.

If $P$ is a principal $G$-bundle then $[p, m]$ denotes the point in the
associated bundle $P\times _G M = (P\times M)/G$ which is the orbit of
$(p,m) \in P\times M$.

\section{What is a contact fiber bundle?}

Let $F$ be a manifold with a contact distribution $\xi^F \subset TF$.
Since we assume that all contact structures are co-oriented, there is
a 1-form $\alpha^F \in \Omega ^1 (F)$ with $\ker \alpha^F = \xi^F$
which gives $\xi^F$ its co-orientation.  Denote the group of
co-orientation preserving contactomorphisms of $(F, \xi^F)$ by
$\Diff_+ (F, \xi^F)$.  That is, 
$$
\Diff _+ (F, \xi ^F) 
= \{ \varphi : F \to F \mid \varphi^* \alpha^F  =
 e^h \alpha ^F \text{ for some } h\in C^\infty (F)\}.  
$$
\begin{definition}[provisional]\label{def1}
Let $(F, \xi^F)$ be a contact manifold.  A fiber bundle $F \to M
\stackrel{\pi}{\to} B$ is {\bf contact} if the transition maps take
values in $\Diff_+ (F, \xi^F)$.  That is, suppose $\{ U_i\}$ is a
cover of $B$ by sufficiently small sets.  Then there exist
trivializations $\phi_i :\pi\inv (U_i) \to U_i \times F$ such that for
all indices $i,j$ with $U_i \cap U_j \not = \emptyset$ and for every
$b\in U_i \cap U_j$ the diffeomorphisms $\left. \phi_j \circ \phi_i
\inv \right |_{\{b\} \times F} : F \to F$ are elements of $\Diff_+
(F, \xi^F)$.
\end{definition}

\begin{remark}\label{rmrk1}
It follows from the definition that for every point $b\in B$ the fiber
$F_b := \pi \inv (b)$ has a well-defined (co-oriented) contact
structure $\xi^b$.  Namely let $U_i$ be an element of the cover
containing $b$.  Let $pr_2: U_i \times F \to F$ denote the projection
on the second factor. Then $\alpha_i : = (pr_2 \circ \phi_i)^*
\alpha^F \in \Omega (\pi\inv (U_i))$ is a 1-form with the restriction 
$\alpha_i |_{F_b}$ to the fiber being contact. Moreover, if $b$ is
also in $ U_j$ then $\ker (\alpha _i |_{F_b} ) = \ker (\alpha _j
|_{F_b} ) $, so the fiber $F_b$ has a well-defined contact structure $\xi^b$.
We let 
\begin{equation} \label{eq1}
\xi^\nu = \bigcup_{b\in B} \xi^b \subset \cV;
\end{equation}
$\xi^\cV$ is a codimension one subbundle of the vertical bundle $\cV : =
\ker (d\pi : TM \to TB)$.
\end{remark}

\begin{lemma} \label{lem1}
Let $F\to M \stackrel{\pi}{\to} B$ be a contact fiber bundle as
 in Definition~\ref{def1} above.  There exists a one-from
$\alpha$ on $M$ such that $(\ker \alpha) \cap \cV = \xi^\nu$, where
$\cV\subset TM$is the vertical bundle and $\xi^\nu$ is the bundle
defined by (\ref{eq1}).  In other words, a contact fiber bundle has a
globally defined one-form that restricts to a contact form on each
fiber.
\end{lemma}

\begin{proof}
Let $\{U_i\}$ be a sufficiently small open cover of $B$.  Choose a
partition of unity $\rho_i$ subordinate to $\{U_i\}$.  Let $\alpha_i
\in \Omega^1 (\pi\inv (U_i))$ denote the one-forms constructed in
Remark~\ref{rmrk1}.  Then $\alpha = \sum_i (\pi^* \rho_i) \alpha_i$ is
a globally defined one-form on $M$ with $(\ker \alpha )\cap T (F_b)$ a 
contact structure  on each fiber $F_b$.
\end{proof}

Lemma~\ref{lem2} and Corollary~\ref{cor1} below form a converse to
Lemma~\ref{lem1}.  Recall that a {\bf connection } on a fiber bundle
$\pi: M\to B$ is a choice of a complement $\cH$ in $TM$ to the
vertical bundle $\cV$ of $\pi$ so that $TM = \cH \oplus \cV$.

\begin{lemma}\label{lem2}
Let $F\to M \stackrel{\pi}{\to} B$ be a fiber bundle with a
co-oriented codimension 1 distribution $\xi \subset TM$ such that for
each fiber $F_b$ the intersection $\xi \cap T (F_b)$ is a contact
distribution on $F_b$.  Then
\begin{enumerate}
\item there is a natural connection $\cH = \cH (\xi)$ on the fiber bundle 
$F\to M \stackrel{\pi}{\to} B$;
\item the parallel transport with respect to $\cH$ (when it exists) preserves 
the contact structure on the fibers.  Additionally parallel transport is
co-orientation preserving.
\end{enumerate}
\end{lemma}

\begin{proof} Since $\xi$ is co-oriented there is a 1-form $\alpha $ on 
$M$ with $\ker \alpha = \xi$. 
Let $\omega = d\alpha |_\xi$.  Since fiber restrictions $\alpha
|_{F_b}$ are contact, $\omega |_{\xi^\nu}$ is non-degenerate, where as
above $\xi^\nu$ is the intersection of the distribution $\xi$ with the
vertical bundle $\cV$.  We define $\cH$ to be the
$\omega$-perpendicular to $\xi^\nu$ in $(\xi, \omega)$: 
$$
\cH = (\xi^\nu)^\omega .
$$ 
Note that if $\alpha'$ is another 1-form with $\ker \alpha' = \xi$
giving $\xi$ its co-orientation then $\alpha' = e^f \alpha$ for some
function $f \in C^\infty (M)$.  Hence $d\alpha' |_\xi = e^f (d\alpha
|_\xi)$ and consequently the definition of $\cH$ does not depend on
the choice of $\alpha$.  Since $\omega |_{\xi^\nu}$ is non-degenerate,
\begin{equation}\label{eq**}
\xi = \xi^\nu \oplus \cH,
\end{equation}
and, since $\xi^\nu = \xi \cap \cV$, $\cH$ is a connection on $\pi: M\to B$.

We now argue that parallel transport defined by $\cH$ preserves
$\xi^\nu$. (Here we tacitly assume that the parallel transport exists,
i.e., that the connection $\cH$ is complete.  If the fiber $F$ is
compact then $\cH$ is complete, but this need not be true in full
generality.  If the parallel transport doesn't exist globally, it does
exist locally: one can parallel transport for short times small pieces
of the fibers.  The statement of the lemma then becomes messy.  And so
we gloss over this point here and elsewhere in the paper.)  Let $v$ be
a vector field on $B$, let $v^\#$ denote its horizontal lift to $M$
with respect to $\cH$: for each $m\in M$, $v^\# (m)$ is the unique
vector field in $\cH_m \subset T_m M$ with $d\pi (v^\# (m)) = v( \pi
(m))$.  Let $w$ be a section of $\xi^\nu$.  We will argue that the Lie
bracket $[v^\# , w]$ is also a section of $\xi^\nu$.  Since $v^\#$ is a
horizontal lift and $w$ is vertical, the bracket $[v^\#, w]$ is also
vertical: $[v^\#, w] \in \Gamma (\cV)$.  By definition of $\cH$, $0 =
\omega (v^\#, w) = d\alpha (v^\#, w)$.  Since $v^\#, w\in \Gamma
(\xi)$ we have $\iota (w) \alpha = 0 = \iota (v^\#)\alpha$.  Therefore
$0 = d\alpha (v^\#, w) = v^\# (\alpha (w)) - w (\alpha (v^\#)) -
\alpha ([v^\# , w])$.  Hence $\alpha ([v^\# , w]) = 0$, i.e., $[v^\# ,
w]\in \xi$.  Consequently $[v^\# , w] \in \Gamma (\cV ) \cap \Gamma (\xi) =
\Gamma (\cV \cap \xi) = \Gamma (\xi^\nu)$, and so the parallel transport with 
respect to $\cH$ preserves $\xi^\nu$.

Finally we argue that the parallel transport also preserves the
co-orientation.  This is a continuity argument.  Let $\gamma :[0, 1]
\to B$ be a path, $\cP_{\gamma (t)} : F_{\gamma (0)} \to F_{\gamma
(t)}$ be the parallel transport along $\gamma$.  Since $d\cP_{\gamma
(t)} \left (\ker (\alpha |_{F_{\gamma (0)} } )\right) = \ker (\alpha
|_{F_{\gamma (t)}} )$, we have 
$$ 
(\cP_{\gamma (t)})^* (\alpha|_{F_{\gamma (t)} } ) 
= f_t ( \alpha |_{F_{\gamma (0)} } ) 
$$
for some nowhere zero function $f_t \in C^\infty (F_{\gamma (0)})$
depending continuously on $t$.  Since $f_0 = 1 > 0$, $f_t > 0$ for all
$t \in [0, 1]$.
\end{proof}

\begin{definition}
We will refer to the connection $\cH = \cH (\xi)$ defined in 
Lemma~\ref{lem2} above as a {\bf contact connection}.
\end{definition}

\begin{corollary}\label{cor1}
Let $F\to M \stackrel{\pi}{\to} B$ be a fiber bundle with a
co-oriented codimension 1 distribution $\xi \subset TM$ such that for
each fiber $F_b$ the intersection $\xi \cap T (F_b)$ is a contact
distribution on $F_b$ and the associated contact connection is
complete..  Then $F\to M \stackrel{\pi}{\to} B$ is a contact fiber
bundle in the sense of Definition~\ref{def1}.
\end{corollary}

\begin{proof}
Use parallel transport with respect to the contact connection $\cH$ to
define local trivializations of the fiber bundle $\pi:M \to B$.
\end{proof}

\begin{remark}
Corollary~\ref{cor1} shows that we may also {\em define} a {\bf contact
fiber bundle} to be a fiber bundle $F\to M \stackrel{\pi}{\to} B$ with
a co-oriented codimension 1 distribution $\xi \subset TM$ such that $\xi \cap
T(F_b)$ is a contact structure on each fiber $F_b$ of $\pi$.  {\em From now
one we will take this as our definition of a contact fiber bundle as
opposed to Definition~\ref{def1}.}
\end{remark}
We finish this section by pointing out that contact fiber bundles as
defined in the remark above are easily constructible as associated
bundles.

\begin{lemma}\label{lem2.6}
Let $(F, \xi^F)$ be a contact manifold with an action of a Lie group
$G$ preserving the contact structure $\xi^F$ and its
co-orientation. Let $G\to P \to B$ be a principal $G$-bundle.  For any
($G$-invariant) connection $\cH \subset TP$ the distribution
\begin{equation}\label{eq21}
\xi_M := \overline{\cH}\oplus (P\times _G \xi^F)
\end{equation}
makes $M = P\times _G F$ into a contact fiber bundle.  Here
$\overline{\cH}\subset TM$ is the connection on $M\to B$ induced by
$\cH$.
\end{lemma}
\begin{proof}
See proof of Theorem~\ref{lem3.1} below.
\end{proof}

\section{When is a contact fiber bundle a contact manifold?}

We now consider the conditions on a contact fiber bundle $(F\to M
\stackrel{\pi}{\to} B, \xi \subset TM)$ that ensure that $\xi$ is a
contact structure on the total space $M$. We will see that the
question is related to the fatness of the contact connection
$\cH(\xi)$ and the image of the moment map for the action of the
structure group of $\pi: M\to B$ on the fiber $F$.  


Recall that the curvature of a connection $\cH$ on a fiber bundle
$F\to M \stackrel{\pi}{\to} B$ is a two-form $\Curv_\cH$ on $B$ with
values in the vector fields on the fiber: for $x,y \in T_b B$ and
vector fields $v,w \in \chi (B)$ with $v(b) = x$, $w(b) = y$ 
$$
(\Curv_\cH)_b (x, y) := [v^\# , w^\#] - [v,w]^\# \in \chi (F_b), 
$$
where, as in Lemma~\ref{lem2}, $^\#$ denotes the horizontal lift with
respect to $\cH$.

\begin{proposition}\label{prop3.2}
Let $(F\to M \stackrel{\pi}{\to} B, \xi \subset TM)$ be a contact
fiber bundle, $\cH = \cH (\xi)$ the contact connection, and $\xi^\circ
\subset T^*M$ the annihilator of $\xi$.  The distribution $\xi$ is a contact 
structure on $M$ iff for all $m\in M$ and all $0 \not = \eta \in \xi_m
^\circ$
\begin{equation}\label{eq-star}
\langle \eta, [(\Curv_\cH)_b (\cdot, \cdot)] (m) \rangle : 
T_b B\times T_b B \to \R \text{ is nondegenerate},
\end{equation}
where $b = \pi (m)$.
\end{proposition}

\begin{proof}
Choose a 1-form $\alpha$ on $M$ with $\ker \alpha = \xi$.  The
distribution $\xi$ is contact iff $d\alpha |_\xi$ is nondegenerate.
Now $\xi = \cH \oplus \xi^\nu$ (cf.\ (\ref{eq**})) and $d\alpha |_{\xi
^\nu}$ is nondegenerate since $(F\to M \stackrel{\pi}{\to} B, \xi)$ is
a contact fiber bundle.  Consequently $\xi $ is contact iff $d\alpha
|_\cH$ is nondegenerate.  Now fix a point $m\in M$ and two vectors
$x,y \in \cH_m \subset T_m M$.  Choose vector fields $v, w$ on $B$
with $v^\# (m) = x$, $w^\# (m) = y$ ( as above $^\#$ denotes the horizontal
lift).  Then, omitting evaluations at $m$, we have: 
$$ 
d\alpha (x, y) =
d\alpha (v^\#, w^\#) = v^\# (\alpha (w^\#)) - w^\# (\alpha (v^\#)) -
\alpha ([v^\#, w^\#]) = 0 - 0 - \alpha ([v^\#, w^\#] - [v, w]^\#), $$
since $\alpha (u^\#) = 0$ for any vector field $u$ on $B$.  Since
$(\Curv_\cH) (v, w) = [v^\# , w^\#] - [v,w]^\# $, we have 
$$ 
(d \alpha _m |_{\cH_m}) (x, y)=
d\alpha _m (x, y) =
\langle \alpha _m, [(\Curv_\cH)_{\pi (m)} (d\pi (x), d\pi (y))] (m) \rangle.
$$
For any $0 \not = \eta \in \xi^\circ_m $ there is $s\not = 0$ such that 
$\eta = s \alpha _m$.   Therefore $d\alpha _m|_{\cH_m}$ is nondegenerate iff
(\ref{eq-star}) holds for any $0 \not = \eta \in \xi^\circ_m $.
\end{proof}

We now reinterpret (\ref{eq-star}) in terms of contact moment maps and
fatness.

\subsection*{Contact moment maps   }

\noindent
Consider a manifold $F$ with a co-oriented contact structure $\xi^F$.
As mentioned in Remark~\ref{remark1}, the punctured annihilator bundle
$(\xi^F)^\circ \smallsetminus F$ has two components:$(\xi^F)^\circ
\smallsetminus F = (\xi^F)^\circ_+ \sqcup (- (\xi^F)_+^\circ)$.
Consider the Lie algebra of contact vector fields $\chi (F, \xi^F)$ on
$F$.  Recall that the contact vector fields are in one-to-one
correspondence with sections of the line bundle $TF/\xi \to M$: the map
$\chi (F, \xi) \to \Gamma (TF/\xi)$, $v \mapsto v +
\xi$ gives the bijection. This is a standard fact.  See for example 
\cite{McDS}. There is a natural pairing between the points of 
the line bundle $(\xi^F)^\circ \to F$ and the contact vector fields:
\begin{equation}
	\label{eq6-22}
(\xi^F)^\circ \times \chi (F, \xi^F) \to \R, \quad 
((f, \eta), v) \mapsto \langle \eta, v(f)\rangle   
\end{equation}
for all $f\in F$, $\eta \in ((\xi^F)^\circ)_f$ and $v\in \chi (F, \xi^F)$,
where on the right the pairing is between a covector $\eta \in
(\xi^F)_f \subset T^*_f F$ and a vector $v(f) \in T_f F$. 


Suppose a Lie algebra $\fg$ acts on $F$ by contact vector fields,
i.e., suppose there is a representation $\rho: \fg \to \chi (F, \xi)$
(or an anti-representation; this is the usual problem with left
actions and Lie brackets defined in terms of left invariant vector
fields).  Then the moment map $\Psi : \chi (F, \xi)^* \supset (\xi^F)
^\circ \to \fg^*$ should be the transpose of $\rho$ relative to the
pairing (\ref{eq6-22}) and the natural pairing $\fg^* \times \fg \to
\R$.  However for various reasons (see below) in the case of co-oriented 
contact structures it is more convenient to define the contact moment
map for $\rho$ as a map $\Psi : (\xi^F)^\circ_+ \to \fg^*$:
\begin{equation}\label{eq-mm}
\langle  \Psi ( f, \eta ), X \rangle = \langle (f, \eta), \rho (X) \rangle
= \langle \eta , \rho (X) (f)\rangle . 
\end{equation}
Note that $(\xi^F)^\circ_+$ is a symplectic submanifold of the
cotangent bundle $T^*F$.  Suppose the anti-representation $\rho : \fg\to
\chi (F, \xi)$ comes from a (left) action of a Lie group $G$ on $F$ preserving 
$\xi$ and its co-orientation (with Lie algebra of $G$ being $\fg$).
In this case we write $X_F$ for $\rho (X)$.  The moment map $\Psi:
(\xi^F)_+^\circ \to \fg^*$ is simply the restriction to
$(\xi^F)_+^\circ$ of the moment map $\Phi: T^*F \to \fg^*$ for the
lifted action of $G$ on $T^*F$.  The action of $G$ preserves $\xi^F$
and its co-orientation if and only if the lifted action preserves
$(\xi^F)^\circ_+$.  A $G$-invariant 1-form $\alpha^F$ on $F$ with
$\ker \alpha^F = \xi^F$ is a $G$-equivariant section of the bundle
$(\xi^F)^\circ_+ \to F$.  Therefore the composition
\begin{equation} \label{eq34}
\Psi_{\alpha^F} := \Psi \circ \alpha^F : F \to \fg^*
\end{equation}
is a $G$-equivariant map.  We will refer to it as the {\bf
$\alpha^F$-moment map}. Note that by definition
\begin{equation}\label{eq35}
\langle \Psi_{\alpha^F} , X\rangle (f) = \langle \alpha^F _f , X_F (f) \rangle
= (\iota (X_F)\alpha^F) (f)
\end{equation}
for all $X\in \fg$ and all $f\in F$.  This is the ``classical'' definition 
of a contact moment map (cf.\
\cite{Albert}, \cite{Geiges}, \cite{Bany-Mol}). It depends on a
choice of a contact form, unlike $\Psi: (\xi^F)^\circ_+ \to \fg^*$,
which is ``universal.''

\begin{remark} \label{remark3.2}
 The pairing (\ref{eq6-22}) suggests a different way of 
looking at Proposition~\ref{prop3.2}.  Denote by $\xi^b$ the contact
structure on the fiber $F_b$: $\xi^b = \xi \cap T(F_b)$.  Then for a
point $m\in F_b$, a covector $\eta\in (\xi^b)^\circ_m$ and a vector
field $v\in \chi (F_b, \xi^b)$ we have
$$
\langle (m, \eta), v\rangle = \langle \eta , v(m)\rangle \in \R.
$$ Note that the connection $\cH$ allows us to identify $\xi^\circ$
with $\bigcup_b (\xi^b)^\circ$ and consequently $\xi^\circ_+ =
\bigcup_b (\xi^b)^\circ_+$.  Consequently the curvature $\curv_\cH$ of
the connection $\cH$ gives rise to a well-defined skew-symmetric form
$\sigma_\cH$ on the vector bundle $\tilde{\cH} \to \xi^\circ_+$ which
is the pull-back of $\cH \to M$ by the projection $p: \xi^\circ_+ \to
M$.  Namely, since $\cH \simeq \pi^* (TB)$, $\tilde{\cH} = (p\circ
\pi)^* (TB)$.  So given $m\in M$, $\eta \in (\xi^\circ_+)_m$ and $u,
v\in T_b B$, where $b = \pi (m)$,
\begin{equation}
(\sigma_\cH)_{(m, \eta)} (u, v) := \langle (m, \eta), (\curv_\cH)_b (u,
v)\rangle = \langle \eta, [ (\curv_\cH)_b (u, v)] (m)\rangle .
\end{equation}
Thus Proposition~\ref{prop3.2} asserts:
\begin{gather}\label{eq-starr}
\textsf{ The distribution $\xi$ is a contact structure on $M$ if and only if
$\sigma_\cH$ defined above}\\ \textsf{is a symplectic form on the vector bundle
$\tilde{\cH} \to \xi^\circ_+$.}\notag
\end{gather}
\end{remark}

\subsection*{Fatness}

\begin{definition}[Weinstein, \cite{We}]
A connection one-form $A$ on a principal $G$-bundle $G \to P
\stackrel{\pi}{\to} B$ is {\bf fat} at a point $\mu \in \fg^*$ if for
any point $p\in P$ the bilinear map $\mu \circ (\Curv _A )_p: \cH_p^A
\times \cH_p^A \to
\R$ is non-degenerate, where $\Curv_A $ is the curvature of the
connection form $A$ and $\cH_p^A = \ker (A_p : T_p P \to \fg)$ is the
associated horizontal distribution.
\end{definition}

\begin{remark}
If $A$ is fat at $\mu \in \fg^*$, it is fat at every point of the set
$\{ Ad^\dagger (g) (a \mu) \mid g\in G, a>0\}$ (here and elsewhere
$Ad^\dagger$ denotes the coadjoint action).
\end{remark}

\begin{remark}
Fatness is an open condition.  Thus if $A$ is fat at $\mu$, it is fat
at every point of a $G\times \R^+$ invariant neighborhood of $\mu$ in
$\fg^*$.
\end{remark}

\begin{theorem}\label{lem3.1}
Suppose a Lie group $G$ acts [on the left] on a manifold $F$
preserving a contact distribution $\xi^F$ and its co-orientation; let
$\Psi: (\xi^F)^\circ_+ \to \fg^*$ denote the associated moment map.
Let $G\to P {\to} B$ be a principal $G$-bundle.  Given a connection
1-form $A$ on $P$, there exists a co-oriented codimension 1
distribution $\xi $ on the associated bundle $M = P\times _G F \to B$
which intersects the tangent bundle of each fiber $F_b$ in a contact
distribution isomorphic to $\xi^F$.  Explicitly
\begin{equation} \label{eq-xi}
\xi  = \cH \oplus (P\times _G \xi^F),
\end{equation}
where $\cH$ is the connection on $\pi : M \to B$ induced by $A$.

Moreover, the distribution $\xi$ is a {\em contact structure} on $M$
if and only if the connection $A$ is fat at the points of the points
of the image $\Psi ((\xi^F)^\circ _+)$.
\end{theorem}

\begin{proof}
Since the action of $G$ on $F$ preserves $\xi^F$ and its co-orientation, 
$$
\xi^\cV := P \times _G \xi^F
$$ 
is a well-defined co-oriented subbundle of the vertical bundle $\cV
\simeq P\times _G (TF)$ of $M\to B$.  The connection 1-form $A$ defines a
complement $\cH$ to $\cV$ in $TM$.  Therefore the distribution $
\xi $ on $M$ defined by (\ref{eq-xi}) is a co-oriented codimension 1 
distribution.  By construction, for each fiber $F_b$ we have $\xi \cap
T(F_b) = \xi^\cV |_{F_b} \simeq \xi^F$. ( More precisely, for each
point $p\in P$ we have an embedding $\iota _p : F \hookrightarrow M$,
$\iota _p (f) = [p, f]$ (where $[p,f] \in P\times _G F$ denotes the image
of $(p,f) \in P\times F$).  Then $d\iota_p (\xi^F ) = \xi^\cV |_{F_b}$
where $b= \pi (p)$.)

Now suppose that $A$ is fat at the points of $\Psi ((\xi^F)^\circ
_+)$.  By Proposition~\ref{prop3.2} and Remark~\ref{remark3.2} it
enough to show that for any $[p,f, \eta] \in P \times _G
(\xi^F)^\circ_+ = \xi^\circ_+$, the pairing 
$$ 
(\sigma_\cH)_{([p,f,\eta])} = 
\langle [p, f, \eta], 
\left((\curv _\cH)_b (\cdot, \cdot)\right) ([p, f])\rangle :
T_bB \times T_b B \to \R 
$$
is nondegenerate (where $b = \pi ([p, f])$).

The curvature $\curv_A: \cH^A \times \cH^A \to \fg$ of $A$ defines a
2-form $\overline{\curv}_A $ on $B$ with values in the adjoint bundle
$P\times _G \fg$. 
To write out $\overline{\curv}_A \in \Omega ^2 (B, P\times _G \fg)$
explicitly we need a bit of notation.  For a point $b\in B$ and
vectors $x, y \in T_bB$ denote the horizontal lift of $x$ and $y$ to
$\cH^A _p$ by $x^\#_p $ and $y^\# _p$ respectively.  Then 
$$
(\overline{\curv}_A)_b (x, y) = [ p, (\curv_A)_p (x^\#_p , y^\#_p )]
\in P\times _G \fg 
$$ 
for any $p\in P$ in the fiber of $P\to B$ above $b$.

Since $G$ acts on $F$ by contact transformations, there is an
(anti-)representation $\rho : \fg \to \chi (F, \xi^F)$, $\rho (X) =
X_F$, from the Lie algebra $\fg$ of $G$ to contact vector fields on
$F$.  Recall: the moment map $\Psi : (\xi^F)^\circ_+ \to \fg$ is the adjoint
of $\rho$ in the sense that
$$
\langle (f, \eta), \rho(X) \rangle = \langle \Psi (f, \eta), X\rangle 
$$ 
for all $f\in F$, $\eta \in (\xi^F)^\circ _+ $, $X\in \fg$,
where on the right the pairing is the natural pairing $\fg^* \times
\fg \to \R$ and on the left it is the pairing (\ref{eq6-22}).

Since $\rho$ and $\Psi$ are equivariant, they induce maps of associated bundles
$$
\bar{\Psi} : P \times _G (\xi^F)^\circ_+ \to P \times _G \fg^*, \quad
\bar{\Psi} ([p, f, \eta] ) = [p, \Psi (f, \eta)]
$$
and
$$
\bar{\rho} : P \times _G \fg \to P \times _G (\chi (F, \xi^F)), \quad
\bar{\rho} ([p, X]) = [p, \rho (X)].
$$
The two pairings above give rise to fiber-wise pairings:
$$
(P\times _G \fg^* ) \oplus (P\times _G \fg) \to P\times _G \R = B \times \R,
\quad [p, \mu] \oplus [p, X] \mapsto ( \pi (p), \langle \mu, X\rangle)
$$
and
$$
P \times _G ((\xi^F ) ^\circ_+ \times \chi (F, \xi^F)) \to P\times _G \R 
= B \times \R,
\quad [p, ((f, \eta), v)]\mapsto (\pi (p), \langle \eta , v(f)\rangle), 
$$
where $\pi$ now denotes the projection $P \to B$.  The maps
$\bar{\rho}$ and $\bar{\Psi}$ are adjoint with respect to the two
fiber-wise pairings.

Finally, the map $\bar{\rho}$ relates the two curvatures,
$\overline{\curv}_A \in \Omega^2 (P\times_G \fg)$ and $\curv_\cH \in
\Omega ^2 (B, P \times _G \chi (F, \xi^F))$:
$$
\bar{\rho} \circ \overline {\curv}_A = \curv_\cH .
$$

Putting together the above remarks we get 
\begin{align*}
(\sigma_\cH)_{[p,f,\eta]} (u,v) 
        &=& \langle [p,f, \eta ], (\curv_\cH)_b (u, v) \rangle \\
&=& \langle [p,f, \eta ], 
\bar{\rho} \circ (\overline{\curv}_A)_b (u, v) \rangle \\
&= & \langle \bar{\Psi} ([p, f, \eta]), (\overline{\curv}_A)_b (u, v) \rangle \\
&=& \langle \Psi (f, \eta), (\curv_A)_p (u^\# _p, v^\#_p)\rangle
\end{align*}
for all $[p, f, \eta] \in P\times _G (\xi^F )^\circ _+$ and any $u,
v\in T_b B$ ($b = \pi ([p, f]) \in B$).
Thus $A$ is fat at the points of $\Psi ((\xi^F)^\circ_+)$ if and only if 
$\sigma _\cH$ is nondegenerate.
\end{proof}

\begin{remark}
Theorem~\ref{lem3.1} allows us to re-interpret (\ref{eq-starr}).
Namely, suppose $(F \to M \to B, \xi)$ is a contact fiber bundle and
$\cH (\xi)$ is the corresponding contact connection. {\em Suppose the
holonomy group $G$ of $\cH (\xi)$ is a finite dimensional 
Lie group. } 
Then $M$ is an associated bundle for a principal $G$-bundle $G \to P
\to B$ and $\cH(\xi)$ is induced by a connection $A$ on $P$.  Also,
the action of $G$ on a typical fiber $(F, \xi^F)$ is contact and
co-orientation preserving.  Then by Theorem~\ref{lem3.1}, the
distribution $\xi$ is a contact structure if and only if $A$ is fat at
the points of the image of the moment map $\Psi: (\xi^F)^\circ_+ \to
\fg^*$ associated to the action of the holonomy group on a typical
fiber.

In general this gives us a {\em formal} interpretation of (\ref{eq-starr})
as fatness of the connection on the principal $G$-bundle $P\to B$
where $G$ is the group of co-orientation preserving contactomorphisms
$\Diff _+ (F, \xi^F)$ and $P$ is the ``frame bundle'' of the fiber bundle
$M\to B$.
\end{remark}

\begin{remark} \label{rem30}
Suppose $F$ is a manifold with an action of a Lie group $G$ and
$\alpha _F$ is a $G$-invariant 1-form on $F$.  Then a choice of a
connection 1-form $A$ on a principal $G$-bundle $G\to P \to B$ defines
a 1-form $\alpha _M = \alpha _M (A, \alpha_F) $ on the associated bundle
$M: = P\times _G F$ such that $\alpha _M$ restricted to each fiber of
$M\to B$ is $\alpha_F$:

Define a moment map $\Psi_{\alpha_F}: F \to \fg^*$ by
$$
\langle \Psi_{\alpha_F} (f), X\rangle = \alpha _F (X_F)
$$
for all $X\in \fg$, where $X_F$ denotes the vector field induced by
$X$ on $F$ (cf.\ (\ref{eq35})).  It is easy to check that the 1-form
$\alpha$ on $P\times F$ given by 
$$
\alpha _{(p,f )} = \langle \Psi_{\alpha_F} (f) , A_p\rangle + (\alpha _F)_f
$$ 
is basic relative to the diagonal action of $G$ on $P\times F$: 
$$
g\cdot (p, f) = (p\cdot g\inv, g\cdot f).  
$$ 
(Since $P$ is a principal $G$-bundle, the natural action of $G$ on $P$
is a right action: $(g, p) \mapsto p\cdot g$. The diagonal action of
$G$ above is a left action.  This matters because of the signs below.)
Now, for any $X\in\fg$,
$$
\iota (X_{P\times F}) \alpha = \langle \Psi_{\alpha_F}, A (X_P) \rangle + 
\iota (X_F) \alpha _F = 
\langle \Psi_{\alpha_F} , -X\rangle + \langle \Psi_{\alpha_F}, X \rangle = 0.
$$ 
Since $\alpha$ is $G$-invariant, it descends to 1-form $\alpha_M$
on $M$.

By Theorem~\ref{lem3.1}, if $\alpha_F$ is a contact form and if the
connection $A$ is fat at the points of $\Psi_{\alpha_F} (F)$, then $\alpha_M$ is a
contact form. Moreover in this case $\xi = \ker \alpha_M$ is precisely
the distribution on $M$ defined by equation (\ref{eq21}).
\end{remark}

In the rest of the paper we discuss two applications of Theorem~\ref{lem3.1}
--- K-contact fiber bundles and contact cross-sections.

\section{Application 1: $K$-contact fiber bundles}

\begin{definition}
Let $(F, \xi)$ be a contact co-oriented manifold. It is {\bf
$K$-contact} if there is a metric $g$ on $F$ such that
\begin{enumerate}
\item the unit normal $n$ to the contact distribution $\xi$, which is 
defined by the metric $g$ and the co-orientation $\xi^\circ _+ \subset
\xi^\circ$ of $\xi$, is Killing, i.e., $L_n g = 0$;
\item the contact form $\alpha_g$ given by $\alpha_g = g(n, \cdot)$ is 
{\bf compatible} with $g$ in the sense that 
\begin{equation} \label{eq-compat}
d\alpha_g |_\xi = (g|_\xi) (\cdot, J \cdot)
\end{equation}
for some complex structure $J$ on $\xi$ with $J^* (g|_\xi) = g|_\xi$.
\end{enumerate}
We will refer to the triple $(F, \xi^\circ_+ , g)$ as a {\bf $K$-contact
structure} and to $g$ as a {\bf $K$-contact metric} on $(F, \xi)$.
\end{definition}
Note that the vector field $n$ in the definition above is the Reeb
vector field of the contact form $\alpha_g$.  

\begin{remark}
Given a contact form $\alpha$ on a manifold $F$ we can easily find a
metric $g$ on $F$ such that the Reeb vector field $R_\alpha$ of
$\alpha$ is unit and normal to $\xi = \ker \alpha$ and such that
$\alpha$ and $g$ are compatible ((\ref{eq-compat}) holds).  If
$R_\alpha$ happens to be Killing with respect to $g$ then $g$ is a
$K$-contact metric.
\end{remark}

We now relate, following Yamazaki, $K$-contact structures on compact
manifolds and contact torus actions (cf.\ \cite{Yamazaki},
Proposition~2.1).
\begin{proposition}
A compact  contact co-orientable  manifold $(F, \xi)$ 
 admits 
a $K$-contact metric $g$ if and only if there is an action of a torus
$T$ on $F$ preserving $\xi$ and a  vector $X$ in the Lie
algebra $\ft $ of $T$ so that the function $\langle \Psi, X \rangle :
\xi^\circ_+ \to \R$ is strictly positive.  Here $\Psi: \xi^\circ _+ \to
\ft^*$ is the moment map associated with the action of $T$ on $(F,
\xi)$.
\end{proposition}

\begin{proof}
Suppose there is an action of a torus $T$ on $(F, \xi)$ and $X\in \ft$
such that $\langle \Psi, X \rangle : \xi^\circ_+ \to \R$ is strictly
positive.  Since the action of $T$ preserves $\xi$, the lifted action
of $T$ on $T^*F$ preserves $\xi^\circ$.  Since $T$ is connected, the
lifted action preserves a component $\xi^\circ_+$ of $\xi
\smallsetminus F$.  It follows that for any 1-form form $\beta $ on
$F$ with $\ker \beta = \xi$, the average $\bar{\beta} $ of $\beta$
over $T$ still satisfies $\ker \beta = \xi$ (if the action of a group
does not preserve the co-orientation of $\xi$, the average of $\beta$
may be zero at some points).  Hence we may assume that there is a
$T$-invariant 1-form $\alpha'$ with $\alpha' (F) \subset \xi^\circ_+$.
Now let 
$$
\alpha = (\langle \Psi \circ \alpha', X\rangle)\inv \alpha'.
$$ 
Then, since $\iota (X_F) \alpha' = \langle \Psi \circ \alpha',
X\rangle$ (cf.\ equations (\ref{eq34}) and (\ref{eq35})), $\iota (X_F)
\alpha = 1$.  Then $TF = \xi \oplus \R X_F $ and the splitting is
$T$-equivariant.  We use the splitting to define the desired metric
$g$.  We declare $\xi $ and $\R X_F$ to be orthogonal and define
$g(X_F, X_F) =1$, so that $X_F$ is a unit normal to $\xi$.  On $\xi$
we choose a $T$-invariant complex structure $J$ compatible with
$d\alpha|_\xi$ and define $g|_\xi (\cdot, \cdot) = d\alpha |_\xi
(\cdot, J\cdot)$.  Then $g$ is $T$-invariant and hence $L_{X_F} g = 0$.
Thus $g$ is a $K$-contact metric on $(F, \xi)$. 

Conversely, if there is a metric $g$ on $F$ making $(F, \xi)$
$K$-contact, the flow $\{\exp (t n) \}$ of the unit normal vector
field $n$ to $\xi$ is a 1-parameter group of isometries $\Diff (F,
g)$.  Since $F$ is compact $\Diff (F, g)$ is a compact Lie group.
Hence the closure $T = \overline{\{\exp (t n) \}}$ is a compact
connected abelian Lie group, i.e., a torus.  Let $X$ be the vector in
the Lie algebra $\ft$ of $T$ with $X_F =n$.  Let $\alpha = g (n,
\cdot)$.  Then $\langle \Psi \circ \alpha , X\rangle = \iota (X_F)
\alpha = g (n, n) =1$.  Hence $\langle \Psi, X\rangle > 0$.
\end{proof}

\begin{theorem}\label{thm-K-cont}
Let $(F, (\xi^F)^\circ_+ , g_F)$ be a compact $K$-contact manifold.
Let $G\subset \Diff (F, g_F)$ be a group of isometries preserving
$(\xi^F)^\circ_+$.  Let $\Psi: (\xi^F)^\circ_+ \to \fg^*$ denote the
associated moment map.  Suppose a principal $G$ bundle $P \to B$ has a
connection 1-form $A$ which is fat at the points of the image $\Psi
((\xi^F)^\circ_+)$.  Then there exists a $K$-contact structure on the
associated bundle $M = P\times _G F$ compatible with the contact form
$\alpha_M = \alpha _M (A, \alpha_{g_F})$ (the form constructed in
remark~\ref{rem30}).  Here $\alpha _{g_F}$ is the
contact form on $F$ defined by $g_F$ and $(\xi^F)^\circ_+$.
\end{theorem}

\begin{proof}  
As we saw above the isometry group $\Diff (F, g_F)$ is a compact Lie
group. Also the flow $\{\exp (t n) \}$ of the unit normal $n$ is a
subgroup of the isometry group whose closure $T = \overline{\{\exp (t
n) \}}$ is a torus.  Since the normal $n$ is $G$-invariant, $T$ and
$G$ commute inside $\Diff (F, g_F)$.  Therefore the torus $T$ acts
naturally on $M = P\times _G F$: $a\cdot [p, f] = [p, a\cdot f]$ for
all $(p, f) \in P\times F$ and all $a\in T$.

The Reeb vector field $R$ of $\alpha_M$ is tangent to fibers of $M\to
B$ , hence $R |_{F_b}$ is the Reeb vector field of $\alpha_M |_{F_b}$
for any fiber $F_b$.  Consequently $R$ is the vector field induced on
$M$ by the $G$-invariant vector field $n \in \chi (F)^G$.  Hence the
flow of $R$ generates the action of $T$ on $M$.  Therefore the
$K$-contact metric $g$ on $M$ has to be $T$-invariant.  Conversely,
any $T$-invariant metric $g$ on $M$ compatible with $\alpha_M$ is
$K$-contact.  The action of $T$ on $M$ preserves the horizontal
subbundle $\cH \subset TM$ defined by $A$, and it preserves the
symplectic structure $d\alpha_M |_\cH$. Choose a $T$-invariant complex
structure $J_\cH$ on $\cH$ compatible with $d\alpha_M |_\cH$.  The
$g_\cH := d\alpha_M |_\cH (J_\cH \cdot, \cdot)$ is a $T$-invariant
metric on $\cH$.  The $T\times G$-invariant metric $g_F$ on $F$ gives
rise to a $T$-invariant metric $g_\cV$ on the vertical bundle $\cV$ of
$M\to B$.  The metric $g_M := g_\cH \oplus g_\cV$ is a $T$-invariant
metric on $M$ compatible with $\alpha_M$ (recall that $\xi : = \ker
\alpha_M = \cH
\oplus (P\times _G \xi^F)$).  Moreover, the Reeb vector field $R$ of 
$\alpha_M$ is unit, normal to $\xi$ and Killing with respect to $g_M$.
Thus $(M, \xi^\circ_+, g_M)$ is a $K$-contact structure on $F \to M =
P\times _G F \to B$.
\end{proof}

\begin{remark}
There is no natural way to make the map $P\times _T F \to B$ into a
Riemannian submersion relative to the $K$-contact metric on $P\times
_T F$ produced by the Theorem~\ref{thm-K-cont}.  Indeed, if we trace
through the construction of $g_M$ we will see that for any point
$[p,f]$ in the fiber $F_b$ of $P\times _G F\to B$ we have 
\begin{equation}
(g_{\cH})_{[p,f]} (v^\#, u^\#) = \langle \Psi_{\alpha_{g_F} }(f), dA_p
(v^\#, J_\cH u^\#)\rangle 
\end{equation}
for any tangent vectors $u, v \in T_bB$ Here on the left hand side
$v^\#$, $u^\#$ denote horizontal lifts to $\cH_{[p,f]}$.  On the right
hand side $v^\#$, $u^\#$ denote horizontal lifts to $\ker A_p \subset
T_pP$.    Thus the horizontal part of the metric $g_M$
depends on the points in the fiber $F_b$!
%
\end{remark}

\begin{example}
 Let $\Sigma$ be a compact Riemann surface and $\omega\in \Omega^2
 (\Sigma)$ an area form which is integral, i.e., $\int _\Sigma \omega
 \in \Z$. Let $S^1 \to P \stackrel{\pi}{\to} \Sigma$ be the
corresponding principal circle bundle with a connection 1-form $A\in \Omega (P,
 \R)^{S^1}$ satisfying $dA = \pi^* \omega$.  Then $A$ is fat on $\R
 \smallsetminus \{0\}$.  Let $\alpha$ be a contact form on a manifold
 $F$ such that the flow of the Reeb vector field $R_\alpha$ is
 periodic.  For example we may take $F$ to be the odd dimensional
 sphere $S^{2n-1} = \{z\in \C^n \mid ||z||^2 = 1\}$ with the standard
 contact form $\alpha = \sqrt{-1} (\sum z_j d \bar{z}_j - \bar{z}_j d
 z_j)|_{S^{2n-1}}$.  Or we can take $F$ to be the co-sphere bundle $S
 (T^* S^k)$ of a sphere with the contact form defined by the standard
 round metric on $S^k$.  Then the associated bundle $P\times _{S^1} F$
 is a $K$-contact manifold.
\end{example}
The next example is a slight generalization.  It produces $K$-contact
manifolds first constructed by Yamazaki by a ``fiber join'' \cite{Yamazaki}.
\begin{example}
For an $n$-tuple $a= (a_1, a_2, \ldots, a_n)$, $a_j >0$, the ellipsoid
$E_a := \{z\in \C^n \mid \sum a_j |z_j|^2 = 1\} \simeq S^{2n-1}$ is a
hypersurface of contact type in $\C^n$. The corresponding contact form
$\alpha_a $ is given by $\alpha _a:= \sqrt{-1} (\sum z_j d \bar{z}_j
-\bar{z}_j d z_j)|_{E_a}$.  For a generic $a$ the Reeb vector field of
$\alpha_a$ generates the action of the $n$-torus $\bbT^n$.  The image
of $E_a$ under the $\alpha_a$-moment map is the simplex 
$$
\Delta_a = \{ (t_1,\ldots, t_n) \in \R^n \simeq Lie(\bbT^n)^* \mid 
\sum a_j t_j = 1, \quad t_j\geq 0 \}.
$$ 
Suppose $\omega_1$, \ldots, $\omega_n$ are integral symplectic
forms on a compact Riemann surface $\Sigma$ such that $\sum t_j
\omega_j$ is nondegenerate for all $t = (t_1, \ldots, t_n) \in
\Delta_a$.  For example we may pick one integral area form $\omega$ and let 
$\omega _j = \omega$ for all $j$.  Then the principal $\bbT^n$ bundle
$P$ over $\Sigma$ defined by $\omega_1, \ldots, \omega_n$ has a
connection 1-form $A = (A_1, \ldots, A_n) \in \Omega (P,
\R^n)^{\bbT^n}$ with $dA = (\pi^* \omega_1, \ldots, \pi^* \omega_n)$.
The connection $A$ is fat at the points of $\Delta_a$.  Therefore
$P\times _{\bbT^n} E_a$ has a $K$-contact structure.  It is an
$S^{2n-1}$-bundle over $\Sigma$.
\end{example}

\section{Application 2: contact cross-sections}

Let $M$ be a manifold with an action of a compact connected Lie group
$G$ preserving a (co-oriented) contact structure $\xi$ on $M$.  Then
there exists a $G$-invariant 1-form $\alpha$ with $\ker \alpha =
\xi$. (Pick any 1-form $\alpha'$ with $\ker \alpha' = \xi$ and average
it over $G$.  Since $\xi$ is $G$-invariant and since $G$ is connected
the averaged form $\alpha$ satisfies $\ker \alpha = \xi$.)  Denote by
$\Psi_\alpha: M \to \fg^*$ the associated $\alpha$-moment map:
$\langle \Psi_\alpha (x), X\rangle = \alpha_x (X_M (x))$ for all $x\in
M$ and all $X\in \fg$; cf.\ (\ref{eq35}).

Since $G$ is compact, for any point $\mu \in \fg^*$ the isotropy Lie
algebra $\fg_\mu$ of $\mu$ has a $G_\mu$-invariant complement $\fm$ in
$\fg$:
\begin{equation}\label{eq5.1}
\fg = \fg_\mu \oplus \fm \quad (G_\mu\text{-equivariant}).
\end{equation}
Moreover we may choose $\fm$ so that $\mu |_\fm = 0$, i.e., $\mu \in
\fm^\circ$. (Pick a $G$-invariant metric on $\fg$ and let 
$\fm = {\fg_\mu} ^\bot$.)  Then a large $\R^{>0}$-invariant open subset
$\cS$ of $\fm^\circ$ is a slice for the coadjoint action of $G$ at
$\mu$.  For example, if $\mu$ is generic, $\fg_\mu$ is a Cartan
subalgebra and $\cS$ is a Weyl chamber (after some
identifications).

We will need (see \cite{SympFib}, p. 37 for a proof):

\begin{lemma}\label{lemma1}
For any $\eta \in \cS$ the pairing
$$
\omega_\eta: \fm \times \fm \to \R, 
\quad (X, Y) \mapsto \langle \eta, [X, Y]\rangle .
$$
is nondegenerate.
\end{lemma}

\begin{theorem}\label{csthm}
Let $(M, \xi = \ker \alpha, \Psi_\alpha : M\to \fg^*)$ be a contact
$G$-manifold as above, $\mu \in \fg^*$ a point, $\fm \subset \fg$ a
subspace defined by (\ref{eq5.1}) with $\mu \in \fm^\circ$ and $\cS
\subset \fm^\circ $ an $\R^{>0}$-invariant slice for the coadjoint action of $G$.   Define 
$$
\cR = \Psi_\alpha \inv (\cS).
$$
Then
\begin{enumerate}
\item  $\cR$ is a contact submanifold of $(M, \xi)$ which is independent of 
the choice of the contact form $\alpha$ used to define it.

\item $G\cdot \cR$ is an open subset of $M$ diffeomorphic to the associated 
bundle $G\times _{G_\mu} \cR$.

\item For any $x\in \cR$
\begin{equation}\label{eq*}
\xi_x = \fm_M (x) \oplus (\xi_x \cap T_x \cR).
\end{equation}
In particular the restriction of the contact structure $\xi$ to
$G\cdot \cR$ is uniquely determined by the $G_\mu$-invariant contact
structure $\xi^\cR := \xi|_\cR \cap T\cR$.
\end{enumerate}
\end{theorem}
\begin{remark}
We will refer to the contact submanifold $(R, \xi^\cR)$ of $(M, \xi)$
as the contact cross-section.
\end{remark}

\begin{proof}[Proof of Theorem~\ref{csthm}. ]
Since $\Psi_\alpha$ is equivariant, the image $d\Psi_\alpha (T_x M)$
contains the tangent space to the coadjoint orbit $G\cdot \Psi_\alpha
(x)$.  Since $\cS$ is a slice, we have
$$
T_\eta \fg^* = T_\eta \cS + T_\eta (G\cdot \eta)
$$
for any $\eta \in \cS$. 
Hence $\Psi _\alpha$ is transverse to $\cS$, and consequently
$$
\cR := \Psi_\alpha \inv (\cS)
$$ is a submanifold.  Also, by equivariance of $\Psi_\alpha$, $\cR$ is
$G_\mu $-invariant.  Since $\cS$ is a slice, $G\cdot \cS (\simeq
G\times _{G_\mu} \cS)$ is open in $\fg^*$.  Hence $G\cdot \cR =
\Psi_\alpha \inv (G\cdot \cS)$ is open in $M$.  Similarly, it's easy to 
see that $G\cdot \cR = G\times _{G_\mu}\cR$.\\

If $\alpha'$ is another $G$-invariant contact form giving $\xi$ its
co-orientation, then $\alpha' = e^f \alpha$ for some function $f\in
C^\infty (M)$.  Consequently $\Psi_{\alpha'} = e^f \Psi_\alpha$.
Since $\cS$ is $\R^{>0}$-invariant, $\Psi_{\alpha'}\inv (\cS) =
\Psi_\alpha \inv (\cS)$.  Thus the cross-section $\cR$ does not depend on 
the choice of the contact form.\\

Note in passing that $\dim \cR = \dim M - \dim G\cdot \mu$, hence odd.
In particular $\cR$ can be contact.\\

We next argue that 
$$
\fm _M  (x) := \{ X_M (x) \mid X\in \fm\}
$$ 
is contained in the contact distribution $\xi_x$ for all $x\in \cR$.
Indeed for any $X\in \fm$ 
$$
\alpha _x (X_M (x)) = \langle \Psi_\alpha (x), X\rangle 
\in \langle \cS, X \rangle \subset \langle \fm^\circ, X\rangle = \{0\},
$$
hence $\fm_M (x) \subset \xi_x$.\\
 
Fix $x\in \cR$.  Since $\fm_M (x) \subset \xi_x$ and since $T_x\cR
\oplus \fm_M (x) = T_xM$, we {\em cannot} have $T_x \cR \subset \xi_x$.  
Therefore $T_x M = T_x \cR + \xi$, and consequently 
$$
\xi^\cR _x := T_x \cR \cap \xi_x
$$ 
is a codimension one subspace of $T_x\cR$.  The rest of the proof
is an argument that $\xi^\cR$ is indeed a contact structure on
$\cR$. In the mean time observe that we have proved (\ref{eq*}). \\

We first argue that the restriction $d\alpha |_{\fm_M (x)}$ is
nondegenerate for all $x\in \cR$. For this we first compute $d\alpha$ on the
tangent space of a $G$-orbit in M.  Let $x\in M$ be a point, $X, Y\in
\fg$ two vectors, $\eta = \Psi_\alpha (x)$.  Then (omitting evaluation
at $x$) we have: $$ d\alpha (X_M, Y_M) = X_M (\alpha (Y_M)) - Y_M
(\alpha (X_M)) - \alpha ([X_M, Y_M]).  $$ Now $$ X_M (\alpha (Y_M)) =
X_M (\langle \Psi_\alpha , Y\rangle ) =
\langle d\Psi_\alpha (X_M), Y\rangle. 
$$
So
$$
X_M (\alpha (Y_M))\, (x) = \langle (d \Psi_\alpha )_x (X_M (x)),  Y \rangle =
\langle X_{\fg^*} (\eta), Y\rangle = - \langle \eta, [X, Y]\rangle,
$$
where the second equality holds by equivariance of the moment map $\Psi_\alpha$
Similarly,
$$
Y_M (\alpha (X_M)) \, (x) = \langle \eta, [X, Y]\rangle .
$$
Since $[X_M, Y_M] = -([X,Y])_M $ (we are dealing with a left action!), 
we have
$$
- \alpha ([X_M, Y_M])\, (x) = \langle \eta, [X, Y]\rangle .
$$
Putting everything together we get 
$$
d\alpha (X_M, Y_M)\, (x) = - \langle \eta, [X, Y]\rangle - 
\langle \eta, [X, Y]\rangle+ \langle \eta, [X, Y]\rangle = 
-\langle \eta, [X, Y]\rangle.
$$
It now follows from Lemma~\ref{lemma1} that for any $x\in \cS$ the restriction 
$d\alpha |_{\fm_M (x)}$ is nondegenerate.\\

We next argue that for any $x\in \cR$ and any $v\in \xi^\cR_x = (T_x \cR) \cap 
\xi_x)$ and any $X\in \fm$ we have $d\alpha _x (X_M (x), v) = 0$.
Let $V$ be a section of $\xi^\cR \to \cR$ with $V(x) = v$.  Then
$d\alpha (V, X_M) = V(\alpha (X_M)) - X_M (\alpha (V)) - \alpha ([V,
X_M])$.  Now $\alpha (V)= 0$ and $\alpha ([V, X_M]) = 0 $ (since $\xi
$ is $G$-invariant $[X_M , V]$ is a section of $\xi$).  Since $V$ is
tangent to $\cR$, we have $(d\Psi_\alpha)_x (V(x)) \in T_{\Psi_\alpha
(x)} \cS = \fm ^\circ$ for all $x\in \cR$.  Therefore $V(\alpha
(X_M))\, (x) = \langle (d\Psi_\alpha)_x (V(x)), X\rangle \in \langle
\fm ^\circ , X\rangle = 0$ for all $x\in \cR$ since $X\in \fm$.
Thus $d\alpha _x (X_M (x), v) = 0$.  Consequently 
$$
\xi_x^\cR \subset \fm_M (x)^{(d\alpha |_{\xi})}.
$$ 
By dimension count the above inclusion is an equality.  

Since
$d\alpha |_{\fm_M (x)}$ is nondegenerate, $d\alpha |_{\xi^\cR_x}$ is
nondegenerate as well for all $x\in \cR$.  Thus $\alpha |_\cR$ is a
contact form, $\xi^\cR = \ker (\alpha|_\cR)$ is a contact structure and 
$\cR$ is a contact submanifold.
\end{proof}

\section*{Acknowledgments}
The work in this paper was partially supported by the Swiss NSF at the
University of Geneva in May 2002.  I am grateful to the University for
its hospitality.

\end{document}